\RequirePackage{fix-cm}

\documentclass[smallextended]{svjour3}       

\smartqed  

\usepackage{graphicx}
\usepackage{epstopdf}

\usepackage{amsmath}
\usepackage{mathrsfs}

\usepackage{algorithm}
\usepackage{algorithmic}

\usepackage{booktabs}
\usepackage{enumitem}

\begin{document}

\title{Application of Progressive Hedging to Var Expansion Planning Under Uncertainty}
\subtitle{ }


\author{Igor Carvalho         \and
        Tiago Andrade         \and
        Joaquim Dias Garcia   \and
        Maria de Lujan Latorre
}


\institute{All authors are with PSR at:\at
              Praia de Botafogo, 370 - Botafogo, Rio de Janeiro - RJ, 22250-040, Brazil \\
              Tel.: +55-21-3906-2100\\
              \email{\{igor,tiago.andrade,joaquim,lujan\}@psr-inc.com}           
                                          \\
}

\date{April 17, 2020}

\maketitle

\begin{abstract}
This paper describes the application of a Progressive Hedging (PH) algorithm to the least-cost var planning under uncertainty. The method PH is a scenario-based decomposition technique for solving stochastic programs, i.e., it decomposes a large scale stochastic problem into s deterministic subproblems and couples the decision from the s subproblems to form a solution for the original stochastic problem. The effectiveness and computational performance of the proposed methodology will be illustrated with var planning studies for the IEEE 24-bus system (5 operating scenarios), the 200-bus Bolivian system (1,152 operating scenarios) and the 1,600-bus Colombian system (180 scenarios).
\end{abstract}

\section{Introduction}\label{Intro}

The fast worldwide insertion of variable renewable energy (VRE) sources such as wind and solar has brought substantial economic and environmental benefits. However, it has also created some technical challenges, for instance, the need for additional generation reserve and/or storage to manage the VRE intermittency. In many systems, VRE sources also disrupted the usual power flow patterns (around which the transmission networks were designed), creating congestions and voltage issues. The solution of these network problems may require substantial investments in new lines, FACTS devices and var sources. Due to the complexity of real networks (size, nonconvexity, discrete variables, uncertainty on load and generator production etc.),  transmission planners usually apply a two-step hierarchical approach.

\begin{figure}[!ht]
\centering
\includegraphics[width=5in]{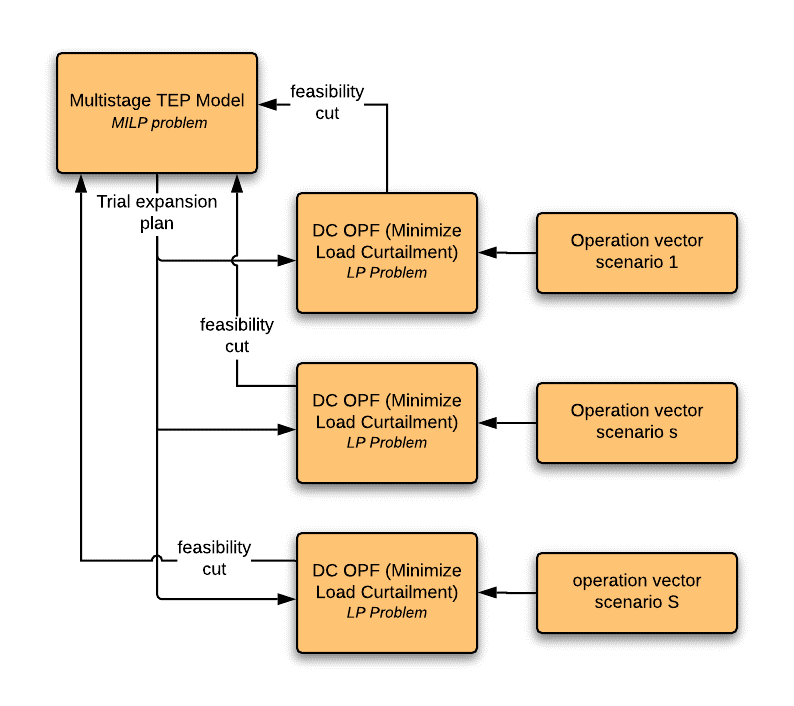}
\caption{First step of the planning process – circuit/transformer reinforcements (Benders decomposition). Source: \cite{latorre2019stochasticrobust}}
\label{fig: first_step}
\end{figure}

The first step determines the least-cost reinforcements of circuits and transformers that eliminate overloads for many operational scenarios (robust optimization). This planning step uses a linearized active optimal power flow (OPF) representation. Because the linearized OPF is convex, it is also possible to use Benders decomposition \cite{latorre2019stochasticrobust}  to handle the multiple operational scenarios (see Fig.~\ref{fig: first_step}).

Given the active power reinforcements from the first step, the second step - which is the focus of this paper - determines the least cost var reinforcements ensuring that voltage levels remain within their limits for the same set of operational scenarios used previously \cite{vaahedi2001dynamic}. This var planning step is arguably more challenging than the first step for two reasons: (i) it requires the solution of AC OPFs, which are nonlinear nonconvex optimization problems; (ii) because of nonconvexities, many decomposition techniques cannot be directly applied.

This paper describes the application of a Progressive Hedging (PH)  algorithm \cite{wets1989aggregation} to the least-cost var planning under uncertainty. Similarly to Benders decomposition \cite{benders1962partitioning}, PH is based on the separate solutions of  AC OPF \cite{granville1988integrated} subproblems for each scenario. The difference is that, instead of an investment module, there is an update in the OPF objective functions. The PH has already been applied successfully to other related problems such as unit commitment \cite{ryan2013toward,gade2016obtaining}, transmission and generation planning \cite{munoz2015scalable}, and a dynamic deterministic OPF \cite{ruppert2018dynamic}, but to the best of our knowledge, it has not been applied yet to the stochastic var expansion problem.

\begin{figure}[!ht]
\centering
\includegraphics[width=5in]{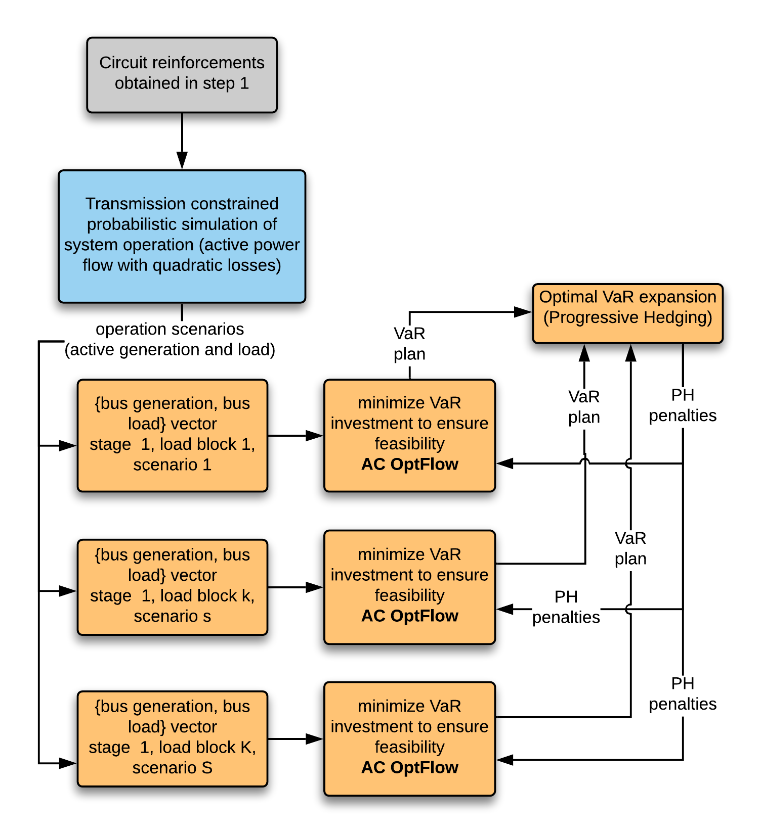}
\caption{PH procedure for var planning under uncertainty.}
\label{fig: PH}
\end{figure}

Although global optimality cannot be guaranteed in the case of nonconvex subproblems, PH has attracted our interest for the following reasons: (i) it is straightforward to implement; (ii) each subproblem is a (slightly modified) OPF, which can be solved by effective techniques such as nonlinear interior point methods \cite{granville1994optimal}; (iii) the OPF subproblems can be solved in parallel; (iv) it has recently become possible to calculate lower bounds to the optimal solution through convex relaxations; and (v) as will be illustrated in this paper, the resulting var plans are near-optimal. Fig.~\ref{fig: PH} shows the proposed PH scheme.

The proposed algorithm is be detailed described in Section~\ref{sec: PH}. Then we present a case study in Section~\ref{sec: study}, where the methodology is initially illustrated for the IEEE 24-bus test case, and then applied to realistic studies of the Bolivian and Colombian power systems. Finally, conclusions are drawn in Section~\ref{sec: conclusions}.

\section{Progressive Hedging Algorithm}
\label{sec: PH}

The PH scheme is implemented in the following steps:

\begin{enumerate}[label=(\alph*)]
	\setlength\itemsep{1em}
	\item Initialize the PH linear and quadratic terms
	
	\item Solve the modified OPF problems for scenarios $s \in 1,\dots, S$. For completeness, we highlight the non-linear OPF model:
	\label{step: Ph_b}
	
	\vspace{\baselineskip} 
	\textbf{Objective function}:
	\begin{align}
    \min & \sum_k [I_k^R Q_{k,(s)}^{C,INV} + I_k^C  Q_{k,(s)}^{C,INV} \\
    & + w_{k,(s)}^{R} (Q_{k,(s)}^{R,INV} - Q_{k}^{R,AVE}) + w_{k,(s)}^{C} (Q_{k,(s)}^{C,INV} - Q_{k)}^{C,AVE}) \\
    & + \rho_k^{R} (Q_{k,(s)}^{R,INV} - Q_{k}^{R,AVE})^2 /2 + \rho_k^{C} (Q_{k,(s)}^{C,INV} - Q_{k}^{C,AVE})^2 /2]
    \end{align}
	
	where:
	
	$I_k^R$ and $I_k^C$ are, respectively, the unit investment costs (\$/Mvar) of reactor and capacitor banks in bus k;
	
	$Q_{k,(s)}^{R,INV}$ and $Q_{k,(s)}^{C,INV}$ are the investment decisions: reactor and capacitor sizes (Mvar) for each scenario s.;
	
	$w_{k,(s)}^R$ and $w_{k,(s)}^C$ are the weights of the linear PH terms;
	
	$\rho_k^R$ and $\rho_k^R$ are the weights of the quadratic PH terms.
	
	As seen in the equation above, the PH terms are related to the linear and quadratic differences between the var investment decisions for each scenario $s$, $Q_{k,(s)}^{R,INV}$ and $Q_{k,(s)}^{C,INV}$, and  “target” investments $Q_{k)}^{R,AVE}$ and $Q_{k}^{C,AVE}$ (known values).

	\vspace{\baselineskip} 
	\textbf{Capacity limits on var injections}:
	\begin{align}
    0 \leq Q_{k,(s)}^{R,INJ} \leq Q_{k,(s)}^{R,INV} \\
    0 \leq Q_{k,(s)}^{C,INJ} \leq Q_{k,(s)}^{C,INV}
    \end{align}
	
	where $Q_{k,(s)}^{R,INJ}$ and $Q_{k,(s)}^{R,INJ}$ are the operational decisions. They are the decision variables to choose how much power will be injected in each bus $k$ for each scenario $s$. They are limited by how much was invested.

	\vspace{\baselineskip} 	
	\textbf{AC power flow equations and constraints}:
	{\small
	\begin{align}
    P_{k,(s)}^{GEN} - P_{k,(s)}^{DEM} - \sum_{j \in \Omega_k} P_{k,j,(s)}^{FLW} (v_k, v_j, \theta_k, \theta_j, t_{k,j}, \phi_{k,j}) = 0 \\
    Q_{k,(s)}^{GEN} - Q_{k,(s)}^{DEM} - \sum_{j \in \Omega_k} Q_{k,j,(s)}^{FLW} (v_k, v_j, \theta_k, \theta_j, t_{k,j}, \phi_{k,j}) + Q_{k,(s)}^{C,INJ} - Q_{k,(s)}^{R,INJ} = 0
    \end{align}    }
    \begin{align}
    \sqrt{{P_{k,j,(s)}^{FLW}}^2 + {Q_{k,j,(s)}^{FLW}}^2} \leq \bar S_{k,j}
    \end{align}
    
    \textbf{Variable bounds}:
    \begin{align}
    \underline Q_{k,(s)}^{GEN} & \leq Q_{k,(s)}^{GEN} \leq \bar Q_{k,(s)}^{GEN} \\
    \underline v_k & \leq v_k \leq \bar v_k \\
    \underline t_k & \leq t_k \leq \bar t_k \\
    \underline \phi_k & \leq \phi_k \leq \bar \phi_k
    \end{align}

    where:
    
    $P_{k,(s)}^{GEN}$ is the active power generation at bus $k$ (known value, taken from the probabilistic simulation of system operation for scenario $s$, see Fig.~\ref{fig: PH});
    
    $Q_{k,(s)}^{GEN}$, reactive power generation at k (decision variable, with limits $\underline Q_{k,(s)}^{GEN}$ and  $\bar Q_{k,(s)}^{GEN}$);
    
    $P_{k,(s)}^{DEM}$ and $P_{k,(s)}^{DEM}$ active and reactive loads at bus k, also known values, taken from the same simulation;
    
    $P_{k,j,(s)}^{FLW}$ and $Q_{k,j,(s)}^{FLW}$ active and reactive power flows in circuit k-j, with limit $\bar S_{k,j}$;
    
    $\theta_k$ and $\theta_j$ voltage angles at buses $k$ and $j$;
    
    $t_{kj}$, tap position value for transformer $k-j$;
    
    $v_k$, voltage at bus $k$, with limits $\bar v_k$ and $\underline v_k$;
    
    $\phi_{kj}$ phase shifting angle in circuit $k-j$, with limits $\bar \phi_{kj}$ and $\underline \phi_{kj}$;
    
    
    

	\item Let  $ \{(\tilde q_k^{rs},\tilde q_k^{cs}), k \in 1,\dots,K \}$ for $s \in 1,\dots, S$ be the optimal solutions of the modified OPF problems of step \ref{step: Ph_b}; carry out the following updates:
	
	\begin{align}
        Q_{k}^{R,AVE} \leftarrow \sum_s Pr(s) Q_{k,(s)}^{R,INV} \\
        Q_{k}^{C,AVE} \leftarrow \sum_s Pr(s) Q_{k,(s)}^{C,INV}
            \end{align}
    \begin{align}
        w_{k,(s)}^R \leftarrow w_{k,(s)}^R + \delta_{k)}^{R} \rho_k^r (Q_{k,(s)}^{R,INV} - Q_{k}^{R,AVE}) \\
         w_{k,(s)}^C \leftarrow w_{k,(s)}^C + \delta_{k)}^{C} \rho_k^r (Q_{k,(s)}^{C,INV} - Q_{k}^{C,AVE})
     \end{align}

    where $\delta_{k)}^{R}$ and $\delta_{k)}^{C}$ are user-defined step sizes, and $Pr(s)$ is vector of the occurrenc probability of each scenario s. It is important to highlight that this step size is a contribution to this paper since the standard PH use an implicit step fixed at 1. In the present work, all scenarios were considered as equiprobable. If the probabilities of the scenarios were different, the method would work in the same way, where the difference consists of the decision variables tending to the value resulting from the weighted average with different weights. 
    
    \item In case of convergence, stop; otherwise go to \ref{step: Ph_b}.
\end{enumerate}

\subsection{Penalty Values Calculations}

The PH performance is based on the chosen values for $\rho_k^R$ and $\rho_k^C$.  With larger values, convergence is faster; on the other hand, there may be an oscillatory behavior. Conversely, smaller values lead to a steadier, but slower, convergence. In this work, $\rho_k^R$ and $\rho_k^C$ were made proportional to the respective investment costs $I_k^R$ and $I_k^C$.

	\begin{align}
        \rho_k^R = K \cdot I_k^R \\
        \rho_k^C = K \cdot I_k^C
    \end{align}
    
where $K$ is a constant chosen to implement variations on the penalty value for each simulation.

\subsection{Convergence}


Two methods are used to measure convergence in this paper. The first one is the usage of a normalized average per-scenario deviation (NApSD) \cite{watson2011progressive} from the reference value detects the proximity to the standard convergence of the method:

\begin{align}
NApSD_{(i)}^{R} = \frac{1}{s}\sum_{k,s | Q_{k}^{R,AVE} > 0} \frac{| Q_{k,(s),(i)}^{R,INV} - Q_{k,(i-1)}^{R,AVE}  |}{Q_{k,(i-1)}^{R,AVE}} \\
NApSD_{(i)}^{C} =\frac{1}{s}\sum_{k,s | Q_{k,(i-1)}^{C,AVE} > 0} \frac{| Q_{k,(s),(i)}^{C,INV} - Q_{k,(i-1)}^{C,AVE}  |}{Q_{k,(i-1)}^{C,AVE}}
\end{align}

The NApSD is calculated for each decision variable and varies at each iteration. Note that the index $i$ indicates the iteration and $i-1$ the previously iteration. As shown above, the NApSD is the average of the relative differences between the investment for scenario $s$ and the average of the investments from the previous iteration. The second method to detect convergence is a relative duality gap computed with valid upper ($UB$) and lower bounds ($LB$). The duality gap is defined as $\frac{UB - LB}{UB}$.

The upper bound is constructed by the union of the investments of all scenarios since an over investment does not harm feasibility. Moreover, Gade et al. \cite{gade2016obtaining} showed how to obtain a lower bound using the PH that is equivalent to the bound provided by the Lagrangian relaxation of the problem that we use in this paper. Since the problem is nonconvex, a zero duality gap could be assured only in special cases \cite{lavaei2011zero}. The average of these problems is equivalent to a Lagrangian relaxation of the full space problem, thus, generating a lower bound (LB). There are alternative methods to obtain a LB, such as \cite{boland2018combining}. The method described in \cite{gade2016obtaining} was implemented since it is simpler. The LB is obtained by solving the PH problem replacing the objective function by the following expression and taking the average for all scenarios.

	\begin{align}
    \min & \sum_k [I_k^R Q_{k,(s)}^{C,INV} + I_k^C  Q_{k,(s)}^{C,INV} \\
    & + w_{k,(s)}^{R} (Q_{k,(s)}^{R,INV} - Q_{k}^{R,AVE}) + w_{k,(s)}^{C} (Q_{k,(s)}^{C,INV} - Q_{k}^{C,AVE})]
    \end{align}

Two additional stop criteria were also considered. The first one consists on verifying the OPF’s convergence for all scenarios under analysis. If any OPF has not converged, the iterative process must stop. The model checks information from the solver if each one of the OPFs has converged or not at each iteration of the method. The other one is for exceeding the maximum number of iterations as established before the method’s execution (timeout). The analysis of simulations with the system is required before defining a good value for this parameter.

\section{Case study}
\label{sec: study}

The present research considers the application of the methodology described for the IEEE 24-bus test case (10 simulations), and also for the electrical power systems from two southamerican countries: Bolivia (5 simulations) and Colombia (1 simulation).

PSR’s software Optflow \cite{optflowmanual} was used to run the optimal power flow simulations and served as a basis for the PH implementation. The OptFlow model is based on the primal-dual interior point algorithm \cite{granville1994optimal}, which have recognized efficiency to solve non convex problems with a large number of variables and constraints. The computer used to run the 10 simulations for the present work consists of an Intel® Core ™ i7-7700K CPU 4.20 GHz with 64 GB of installed RAM.

Artelys Knitro \cite{nocedal2006knitro} is the local nonlinear solver chosen to obtain the solutions from the OPF previously presented. It solves nonlinear problems by considering integer or continuous variables and is prepared for multiple objective functions and nonlinear constraints. 


Since the continuous problem is already complex to solve due to its nonconvexity, we avoided including integer variables into the model. The required investment is decided by choosing how much Mvar is needed at each bus of the system. After the continuous problem is solved, a post-processing is done to choose which shunt equipment should be installed with predetermined capacities and costs.

The list of candidate buses was previously obtained from a process of shunt allocation. Initially, the model is used to identify candidate buses for reactive support, minimizing the reactive power injections in the system, where we discover which buses are in need of reactive support. Results from simulations with IEEE 24-Bus case illustrated the importance of this step. 

The cost considered for each var equipment came from a report from Ente Operador Regional (EOR)  \cite{eor}. The values considered are illustrated in Table~\ref{tab: var cost}.

\begin{table}[!ht]
	\renewcommand{\arraystretch}{1.3}
	\centering
	\caption{Var sources cost}
	\label{tab: var cost}
	\begin{tabular}{c c}
		\toprule
			Equipament & Cost [kUS\$] \\
			\cline{2-2}
            5 Mvar Capacitor  & 313 \\
            10 Mvar Capacitor & 362 \\
            15 Mvar Capacitor & 418 \\
            15 Mvar Reactor   & 1,810 \\
            60 Mvar Reactor   & 2,171 \\
            Connection Bay    & 3,217 \\
		\bottomrule
	\end{tabular}
\end{table}

\subsection{IEEE 24-Bus System}

\vbox{%
The system consists of:

\begin{itemize}
	\item 24 buses (18 with load, 10 with generating units);
	\item 50 transmission lines;
	\item 7 transformers
\end{itemize}
}

\noindent
The list of parameters considered for simulations with the present system is:

\begin{itemize}
	\item Maximum shunt allocation for each bus: 500 Mvar;
	\item 5 different generation/load set-points considered at each iteration;
	\item Maximum number of iterations: 100;
	\item Convergence gap: 1\%
\end{itemize}

With the IEEE 24-Bus database, the first two simulations were done considering ($\delta_k^{rs}$,$\delta_k^{cs}$) = 0 and  K = 1:

\begin{itemize}
	\item \textbf{Simulation 1:} Reactive investment is allowed at all buses from the system (Full set of shunt candidates);
    \item \textbf{Simulation 2:} Reactive investment is allowed only at the buses that presented reactive investment in iteration 0 (Reduced set of shunt candidates)
\end{itemize}

\begin{figure}[!ht]
\centering
\includegraphics[width=5in]{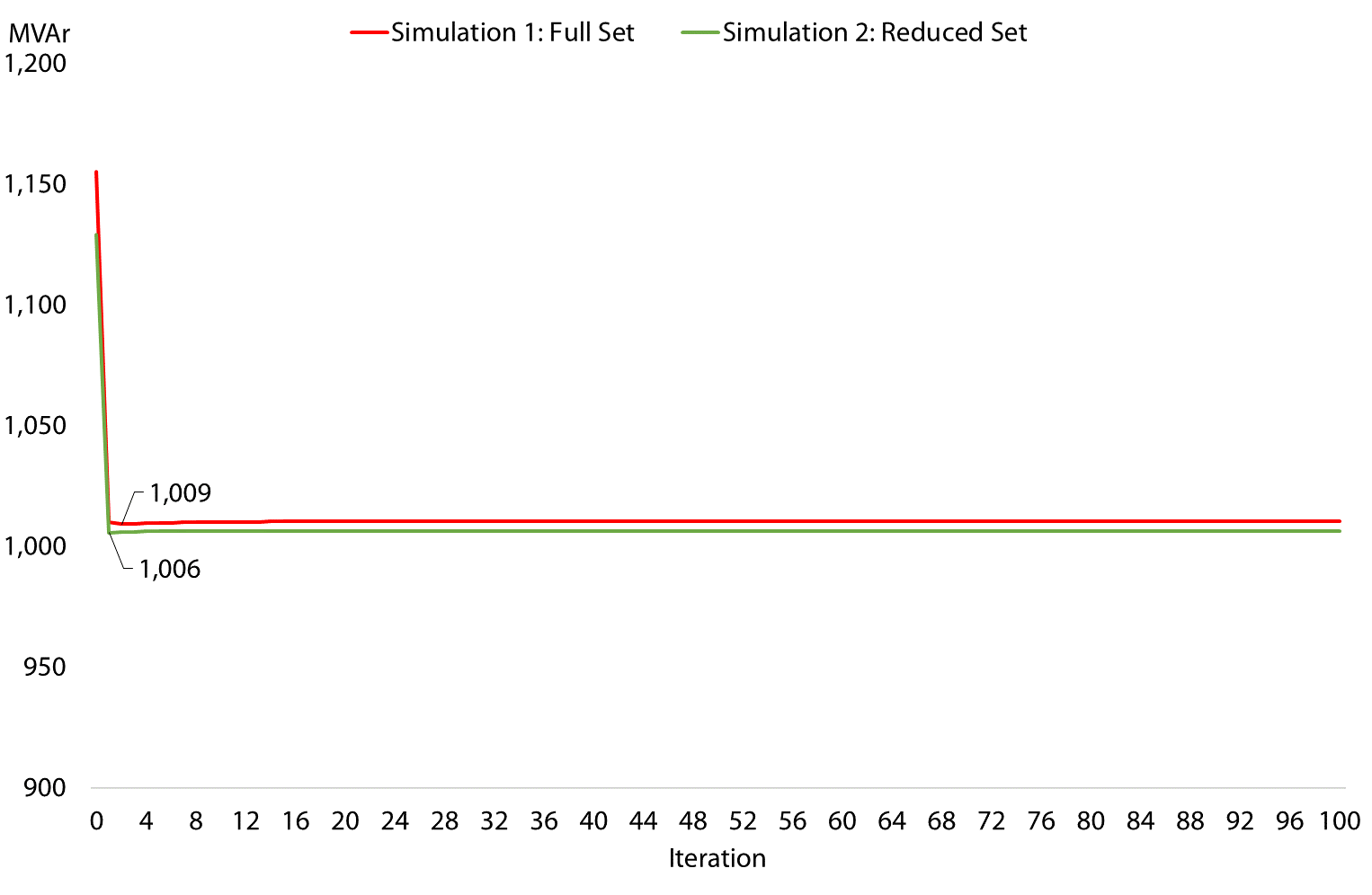}
\caption{Mvar requirement reduction for simulations 1 and 2}
\label{fig: 24-bus Mvar}
\end{figure}

The comparison of the reactive requirement between these two simulations is presented in Fig.~\ref{fig: 24-bus Mvar}.

It is observed in the Fig.~\ref{fig: 24-bus Mvar} that, with the reduced set of candidates, the total amount of shunt requirement reaches slightly smaller values than those from \textbf{Simulation 1}. As the results presented were better for the simulations with a reduced set of candidates, for the next simulations presented in this paper, this reduced set will always be considered.

\vbox{%
In order to evaluate the impact of the step-size, four more simulations were done with K=1. The set to compare now is:

\begin{itemize}
    \item \textbf{Simulation 2:} ($\delta_k^{rs}$,$\delta_k^{cs}$) = 0;
    \item \textbf{Simulation 3:} ($\delta_k^{rs}$,$\delta_k^{cs}$) = 0.001;
    \item \textbf{Simulation 4:} ($\delta_k^{rs}$,$\delta_k^{cs}$) = 0.01;
    \item \textbf{Simulation 5:} ($\delta_k^{rs}$,$\delta_k^{cs}$) = 0.1;
    \item \textbf{Simulation 6:} ($\delta_k^{rs}$,$\delta_k^{cs}$) = 1;
\end{itemize}
}

\begin{figure}[!ht]
\centering
\includegraphics[width=5in]{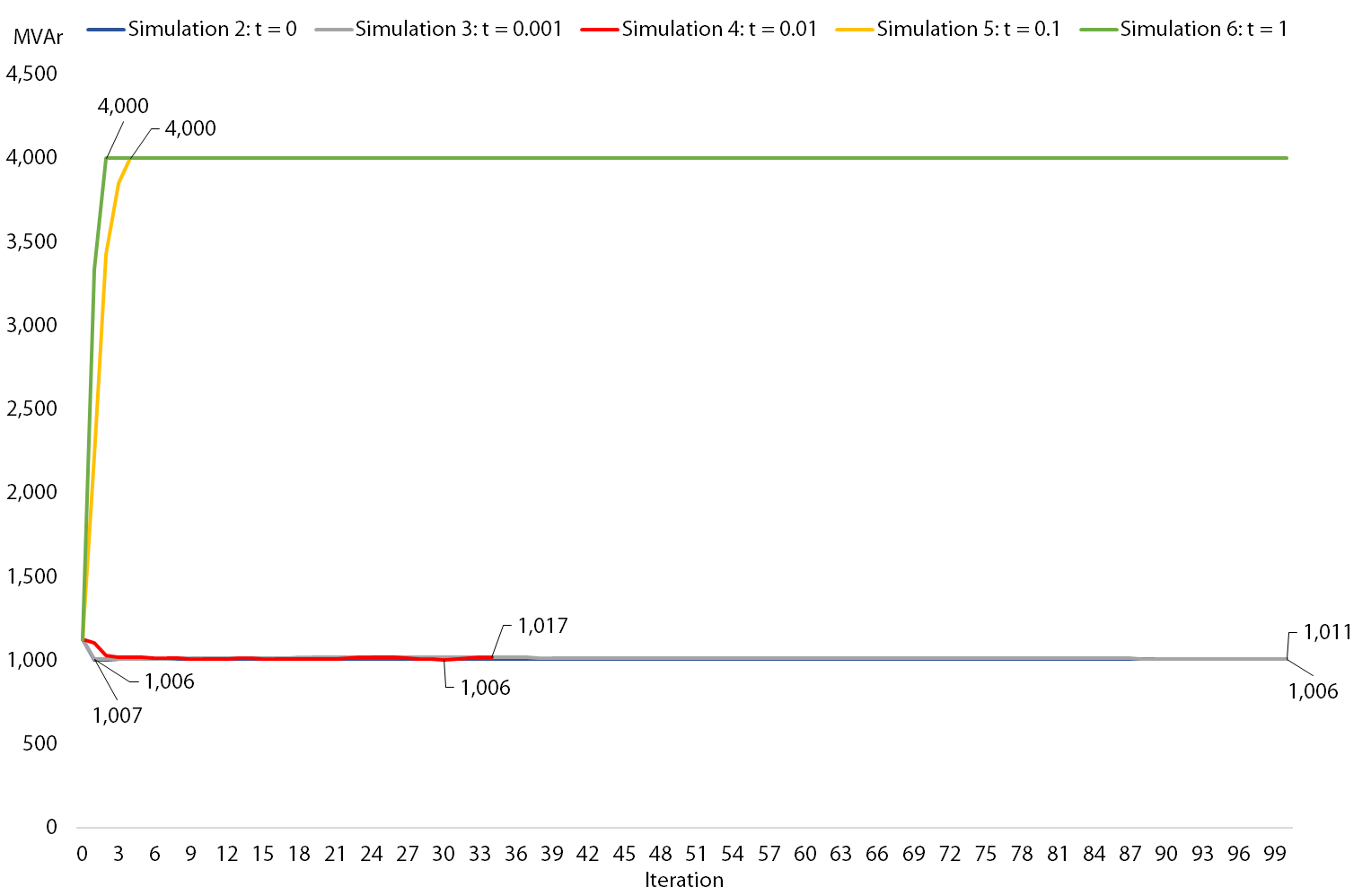}
\caption{Mvar requirement reduction for simulations 2-6}
\label{fig: 24-bus Mvar 2-6}
\end{figure}

The comparison of the total system's reactive requirement along the iterative process for the five simulations above is illustrated in  Fig.~\ref{fig: 24-bus Mvar 2-6}. It is observed that for \textbf{Simulation 5} and \textbf{Simulation 6} (step-sizes 0.1 and 1, respectively), the problem becomes maximizing instead of minimizing in the early first iterations. It occurs because the linear term in the objective function becomes negative in such a way that all objective function becomes negative. Then all reactive investment possible (4,000 MVAr) is allocated. For better comparison between the simulations from this set,  \textbf{Simulation 5} and \textbf{Simulation 6} were both excluded from Fig.~\ref{fig: 24-bus Mvar 2-6}.

\begin{figure}[!ht]
\centering
\includegraphics[width=5in]{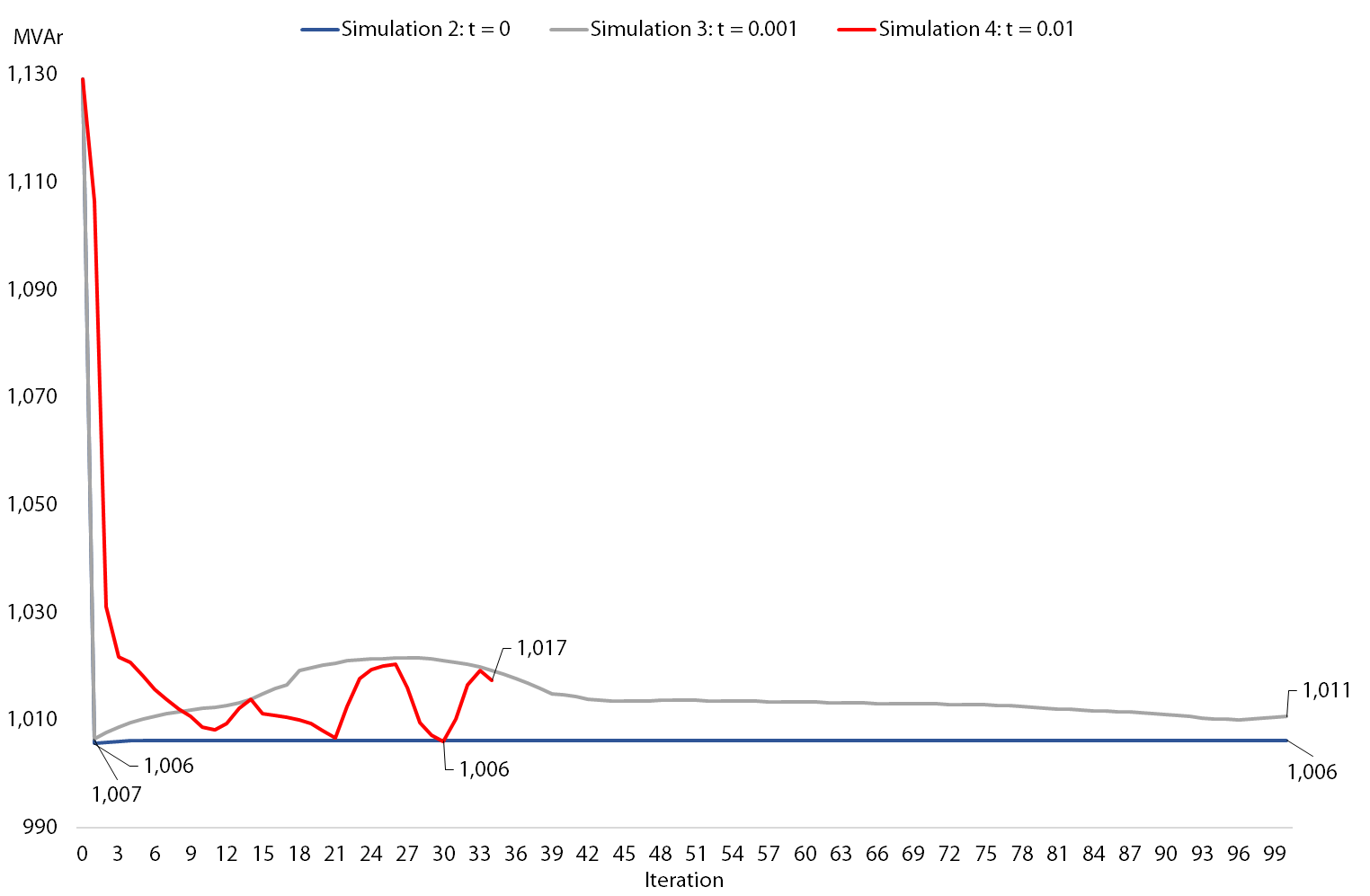}
\caption{Mvar requirement reduction for simulations 2-4}
\label{fig: 24-bus Mvar 2-4}
\end{figure}

\begin{figure}[!ht]
\centering
\includegraphics[width=5in]{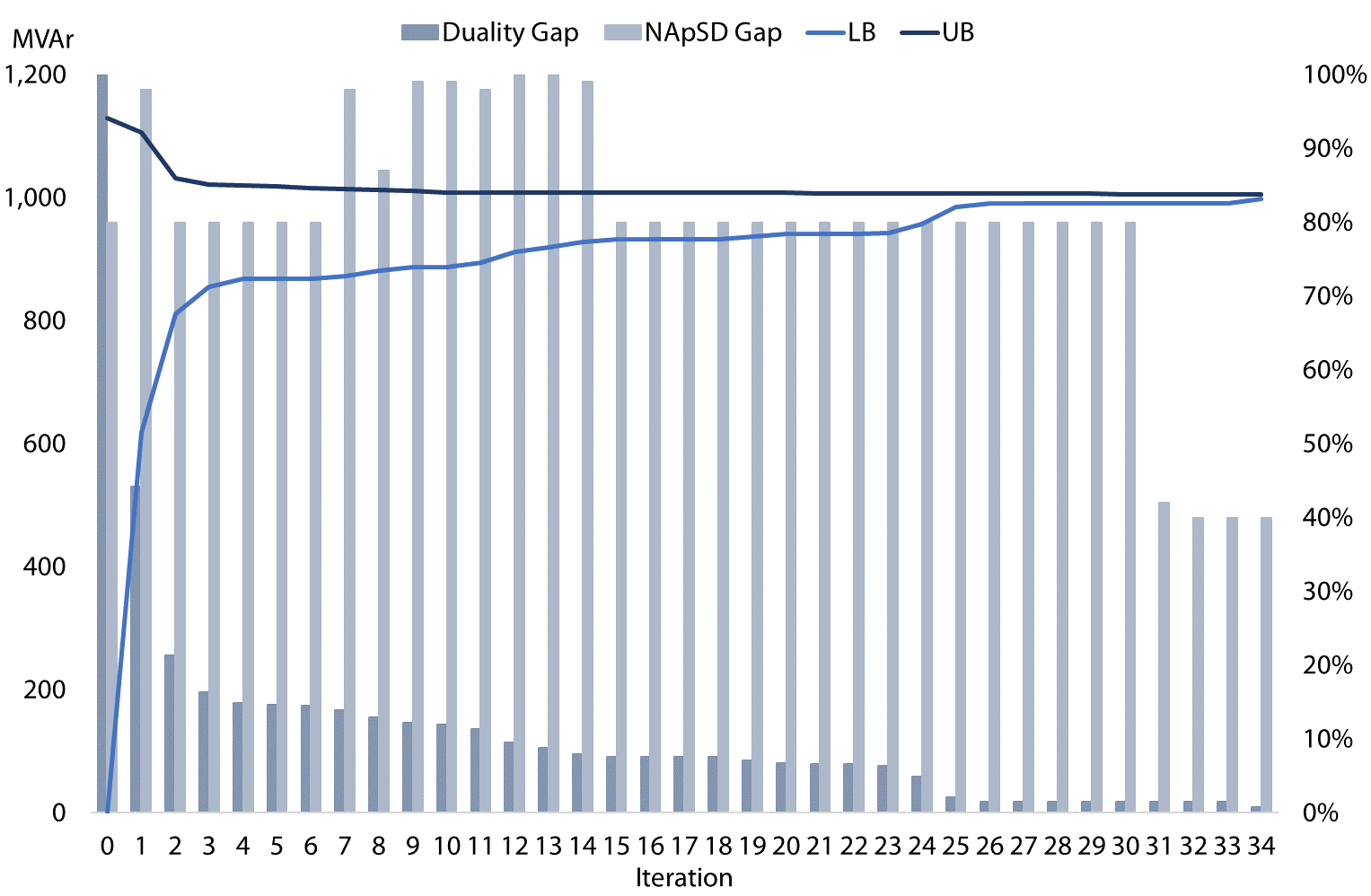}
\caption{LB-UB Convergence scheme – IEE 24-Bus.}
\label{fig: 24-bus LB UB sim4}
\end{figure}

Fig.~\ref{fig: 24-bus Mvar 2-4} illustrates the simulations which reduced the system’s requirement from iteration zero. It is observed that \textbf{Simulation 4} stops before the timeout, because it converges from the duality gap criterium, as illustrated in Fig.~\ref{fig: 24-bus LB UB sim4}.

\begin{figure}[!ht]
\centering
\includegraphics[width=5in]{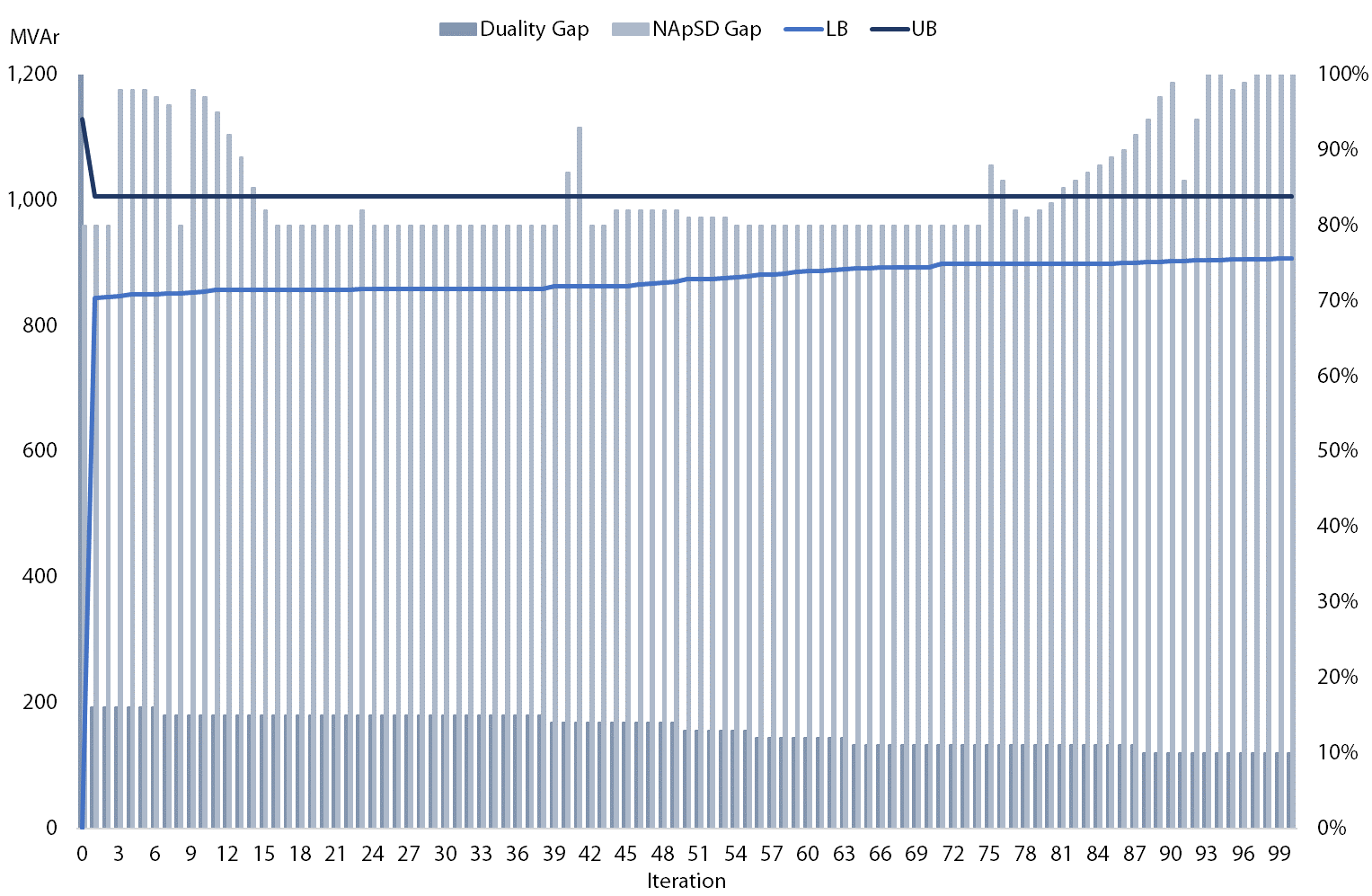}
\caption{LB-UB Convergence scheme – IEE 24-Bus.}
\label{fig: 24-bus LB UB sim3}
\end{figure}

Note that using duality criterion convergence is reached (1\% duality gap) while the NApSD criterion is 40\% as illustrated in Fig.~\ref{fig: 24-bus LB UB sim4}. \textbf{Simulation 4} converges from the duality gap, and \textbf{Simulation 2} does not converge from NApSD. As step-size is null at \textbf{Simulation 2}, it is not possible to compute lower bounds, and therefore, the duality gap is not applicable for this simulation. About \textbf{Simulation 3}, the convergence is slower, as illustrated in Fig.~\ref{fig: 24-bus LB UB sim3}, and it is concluded that this step-size is, on the other hand, too small.

From the previous analyses, it is concluded that the step-size equals to 0.01 is the best option to proceed with further simulations and analyses. However, to reinforce the idea, the same analyses with the step size parameter will be repeated with the Bolivian system in the next section.

\begin{figure}[!ht]
\centering
\includegraphics[width=5in]{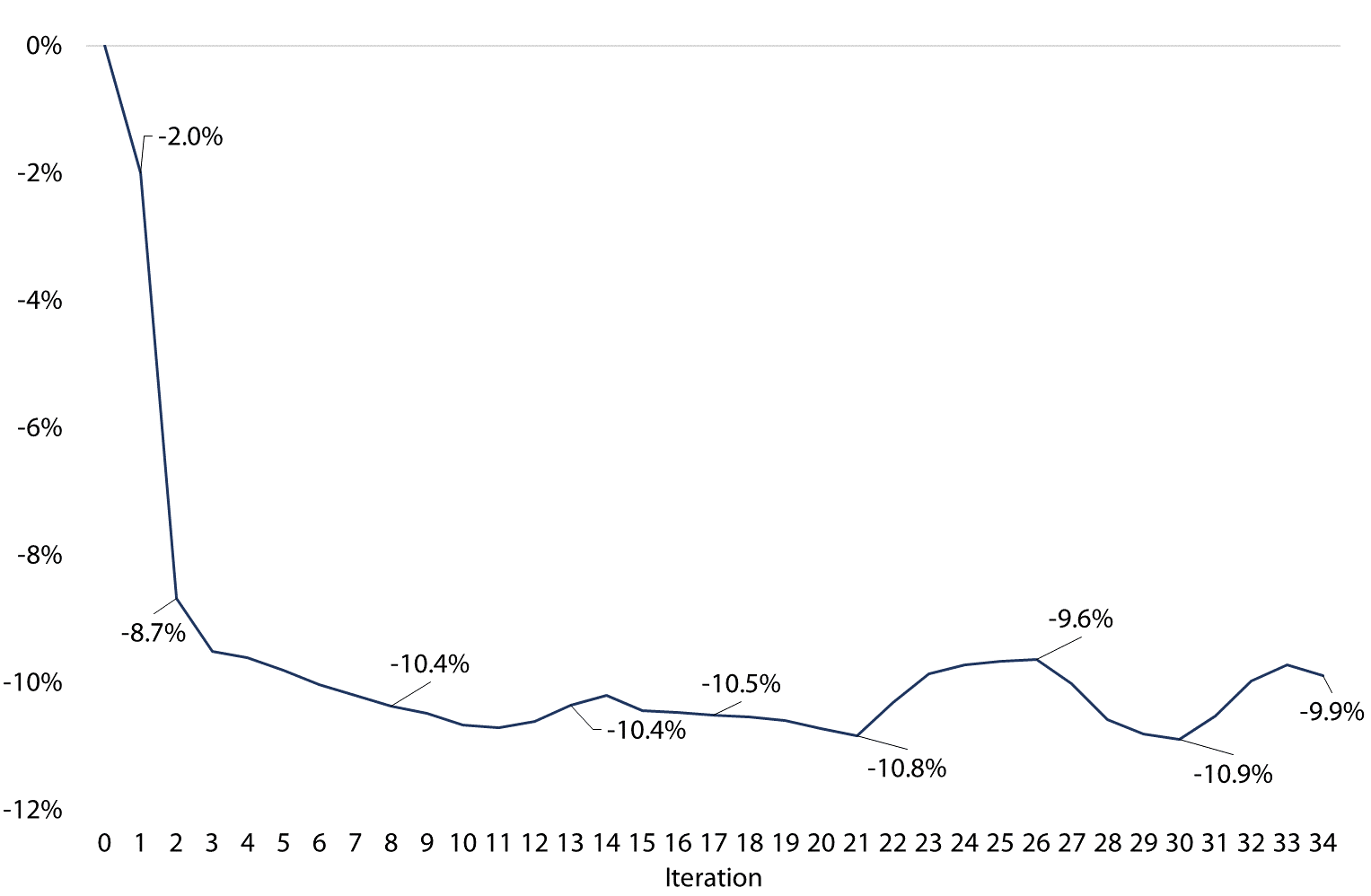}
\caption{Reduction in objective function, \textbf{Simulation 4}– IEEE 24-Bus}
\label{fig: 24-bus of}
\end{figure}

Fig.~\ref{fig: 24-bus of} highlights the reductions on total amount invested compared to the solution obtained without PH (iteration 0). Some oscillatory behavior is observed from variations with weights at each iteration. It is essential to highlight that all solutions are feasible ones, and at the end of the process, the one with the smallest total system's var requirement is chosen. 

\vbox{%
In order to evaluate the impact of the penalty parameter, four more simulations were done with ($\delta_k^{rs}$,$\delta_k^{cs}$) = 0.01. The set to compare now is:
\begin{itemize}
    \item \textbf{Simulation 4:} K = 1;
    \item \textbf{Simulation 7:} K = 0.5;
    \item \textbf{Simulation 8:} K = 2;
    \item \textbf{Simulation 9:} K = 5;
    \item \textbf{Simulation 10:} K = 10;
\end{itemize}
}

\begin{figure}[!ht]
\centering
\includegraphics[width=5in]{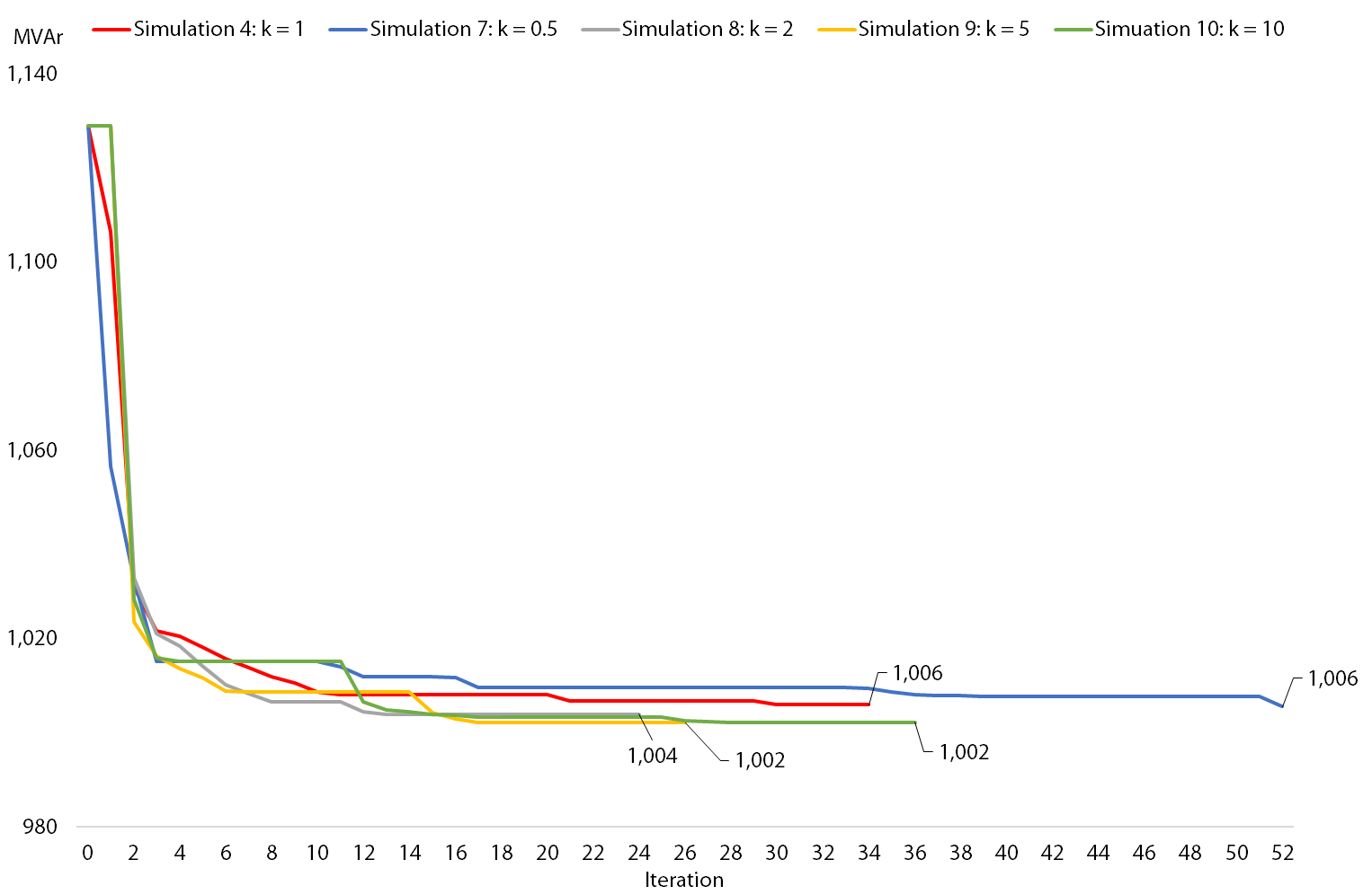}
\caption{Mvar requirement reduction for simulations 4, 7-10}
\label{fig: 24-bus Mvar 4 7-10}
\end{figure}

The comparison of the best total system's reactive requirement found so far along the iterative process for the five simulations above is illustrated in Fig.~\ref{fig: 24-bus Mvar 4 7-10}. It is observed that K = 0.5 or 1 results in similar total system’s var requirement. As the constant K rise from 1 to 5, the local optimal solution found is slightly better. The final solution obtained from simulations 11 and 12 (K = 5 and K = 10, respectively) were also similar.

Table~\ref{tab: CPU Time - IEEE 24} summarizes the total CPU time for the whole iterative process from the 10 simulations analyzed for this system.

\begin{table}[!ht]
	\renewcommand{\arraystretch}{1.3}
	\centering
	\caption{Total CPU Time - IEEE 24 Bus}
	\label{tab: CPU Time - IEEE 24}
	\begin{tabular}{c c c c c c c}
		\toprule
			Simulation & Candidates & t & K & Number of iterations & Total var requirement & CPU Time [min] \\
			\cline{2-7}
            1 & Full set & 0 & 1 & 100 & 1,009 & 4.1 \\
            2 & Reduced set & 0 & 1 & 100 & 1,006 & 4.1 \\
            3 & Reduced set & 0.001 & 1 & 100 & 1,007 & 8.3 \\
            4 & Reduced set & 0.01 & 1 & 34 & 1,006 & 2.8 \\
            5 & Reduced set & 0.1 & 1 & 100 & 1,129 & 8.1 \\
            6 & Reduced set & 1 & 1 & 100 & 1,129 & 8.1 \\
            7 & Reduced set & 0.01 & 0.5 & 52 & 1,006 & 4.3 \\
            8 & Reduced set & 0.01 & 2 & 24 & 1,004 & 2.0 \\
            9 & Reduced set & 0.01 & 5 & 26 & 1,002 & 2.1 \\
            10 & Reduced set & 0.01 & 10 & 36 & 1,002 & 3.0 \\
		\bottomrule
	\end{tabular}
\end{table}

From Table~\ref{tab: CPU Time - IEEE 24} it is observed that considering a full set or a reduced set of candidates makes no visible change in computational effort. It is also observed that with the introduction of the step size, the CPU time after 100 iterations doubles. It happens because another problem is required to be solved to compute lower bounds from the duality gap. Then, with a step size different than zero, two problems are solved at each iteration. With the null step size, only one problem is solved, which explains the rise in the CPU time observed in Table~\ref{tab: CPU Time - IEEE 24}. The simulations considering a value different than zero for the step-size which presented a CPU time much lower than the others, it is because it has converged from the duality gap criterium in the iteration indicated. It is also observable that the value of the penalty parameter seems to have no interference in the CPU Time required to solve the problem with the method’s application.

\subsection{Bolivian Electrical Power System}

\vbox{%
System overview:
\begin{itemize}
	\item 215 buses;
	\item 184 transmission lines;
	\item 59 transformers;
	\item 135 generating units;
	\item 14 shunt equipment
\end{itemize}
}

\noindent
The list of parameters considered for simulations with the present system is:

\begin{itemize}
	\item Maximum shunt allocation for each bus: 500 Mvar;
	\item 1,152 different generation/load set-points considered at each iteration;
	\item Maximum number of iterations: 30;
	\item Convergence gap: 1\%
\end{itemize}

With the Bolivian system database, two sets of simulations were done in order to evaluate the impact of the step-size and the penalty parameter, as previously done with the IEEE 24-bus system. Then, as stated before, all simulations with the Bolivian system considered a reduced set of shunt candidates. First, in order to evaluate the impact of the step-size, three simulations were done with K=1. The set to compare now is:

\begin{itemize}
    \item \textbf{Simulation 11:} ($\delta_k^{rs}$,$\delta_k^{cs}$) = 0.005;
    \item \textbf{Simulation 12:} ($\delta_k^{rs}$,$\delta_k^{cs}$) = 0.01;
    \item \textbf{Simulation 13:} ($\delta_k^{rs}$,$\delta_k^{cs}$) = 0.1;
\end{itemize}

\begin{figure}[!ht]
\centering
\includegraphics[width=5in]{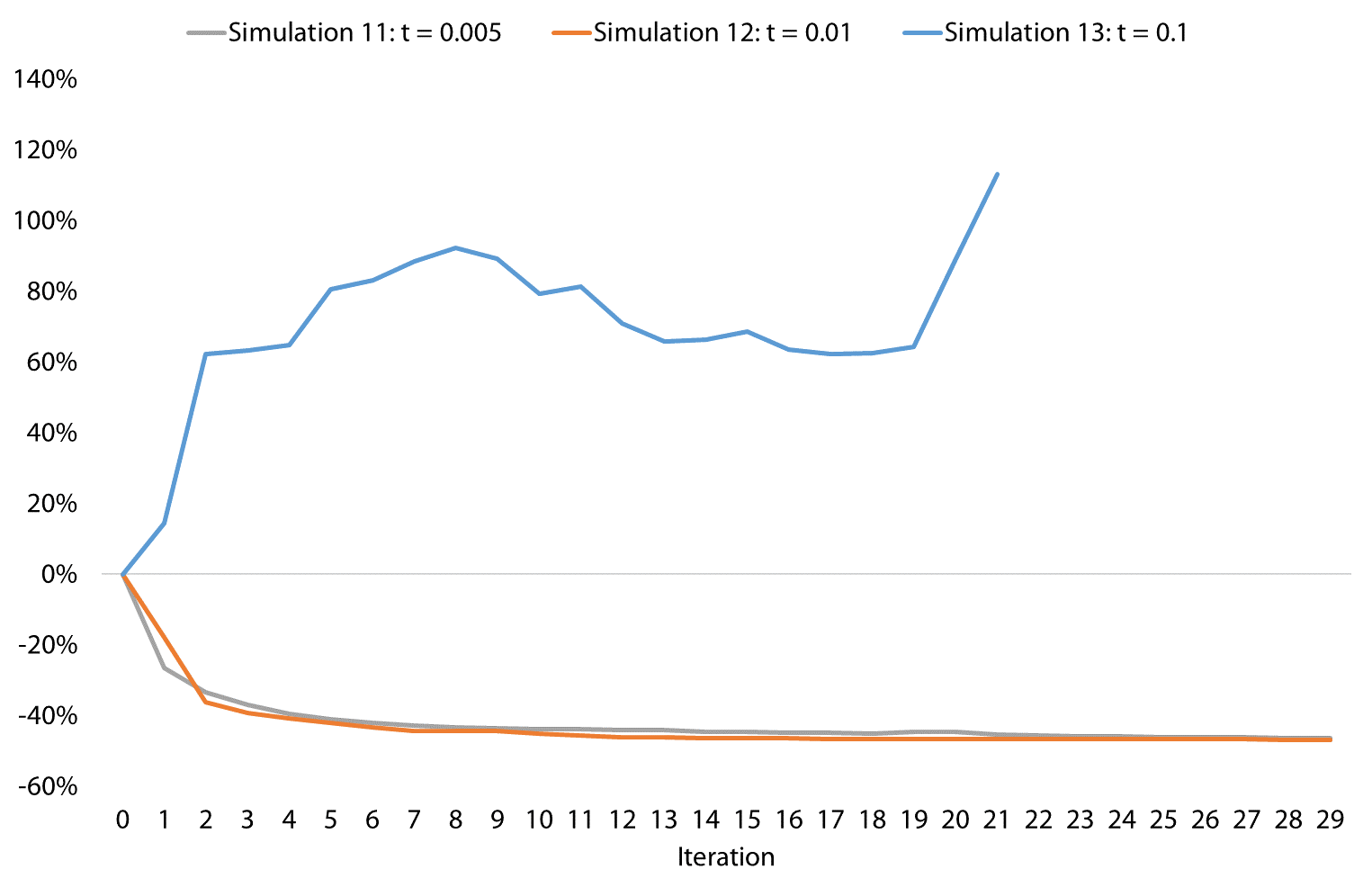}
\caption{Mvar requirement reduction in \% for simulations 11-13 – Bolivia}
\label{fig: bol_stepanalysis}
\end{figure}

The comparison of the reactive requirement between these three simulations is presented in Fig.~\ref{fig: bol_stepanalysis}.

\begin{figure}[!ht]
\centering
\includegraphics[width=5in]{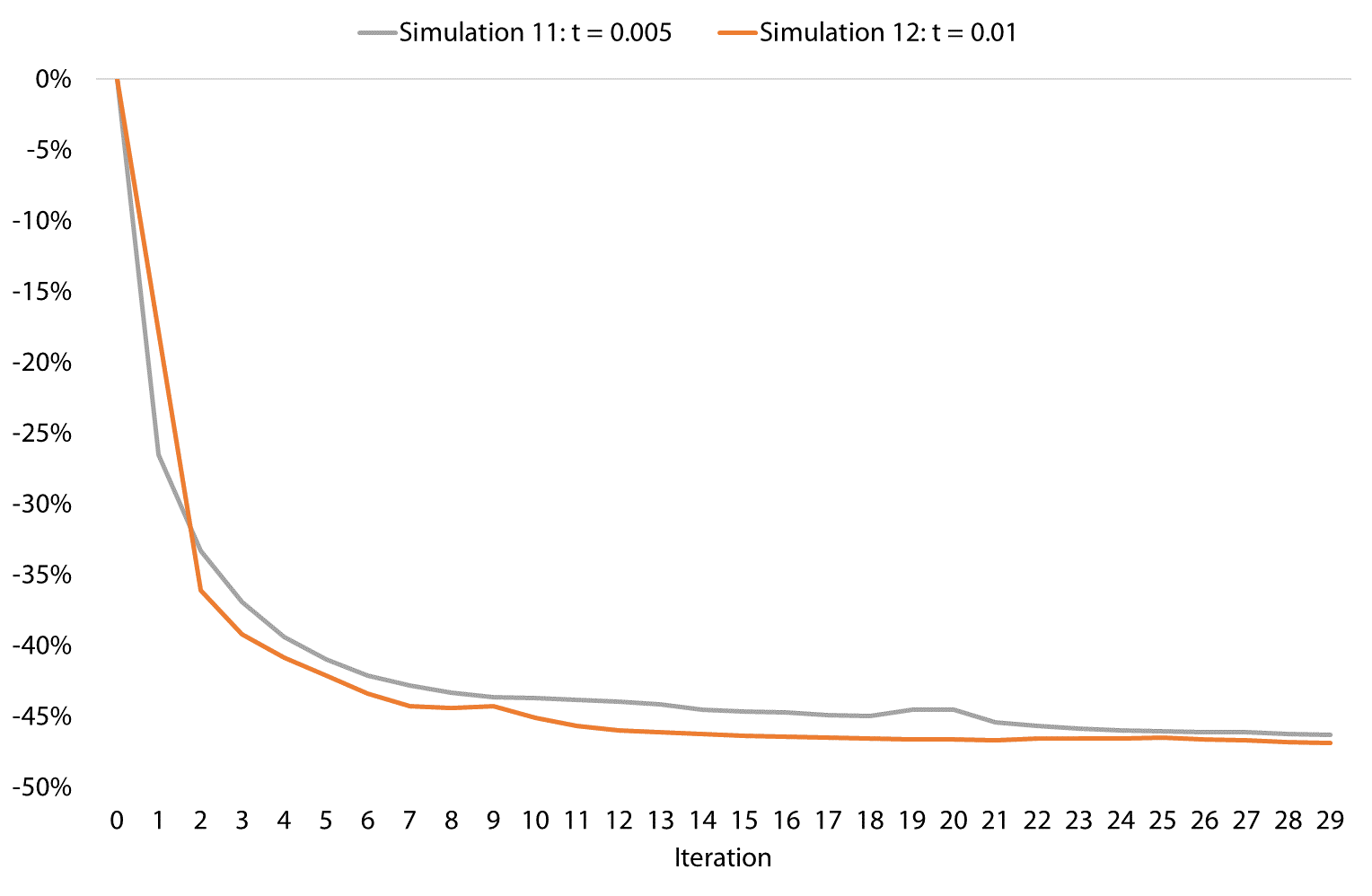}
\caption{Mvar requirement reduction in \% for simulations 11-12 – Bolivia}
\label{fig: bol_stepanalysis2}
\end{figure}

It is observed in Fig.~\ref{fig: bol_stepanalysis} that with a larger step-size (\textbf{Simulation 13}), the weight vectors become negative and higher than the cost in absolute value. It turns the objective function to maximize instead of minimizing the investment in shunt equipment, as observed in the analyses with IEEE 24 bus system. For better comparison, \textbf{Simulation 13} was excluded from the Fig.~\ref{fig: bol_stepanalysis}, and the “simulations with minimization” are illustrated in Fig.~\ref{fig: bol_stepanalysis2}. It is observed that the smaller step-sizes used for \textbf{Simulation 11} and \textbf{Simulation 12} resulted in quite similar amount reductions, where the step-size of 0.01 from \textbf{Simulation 12} presented slightly better results.

\begin{figure}[!ht]
\centering
\includegraphics[width=5in]{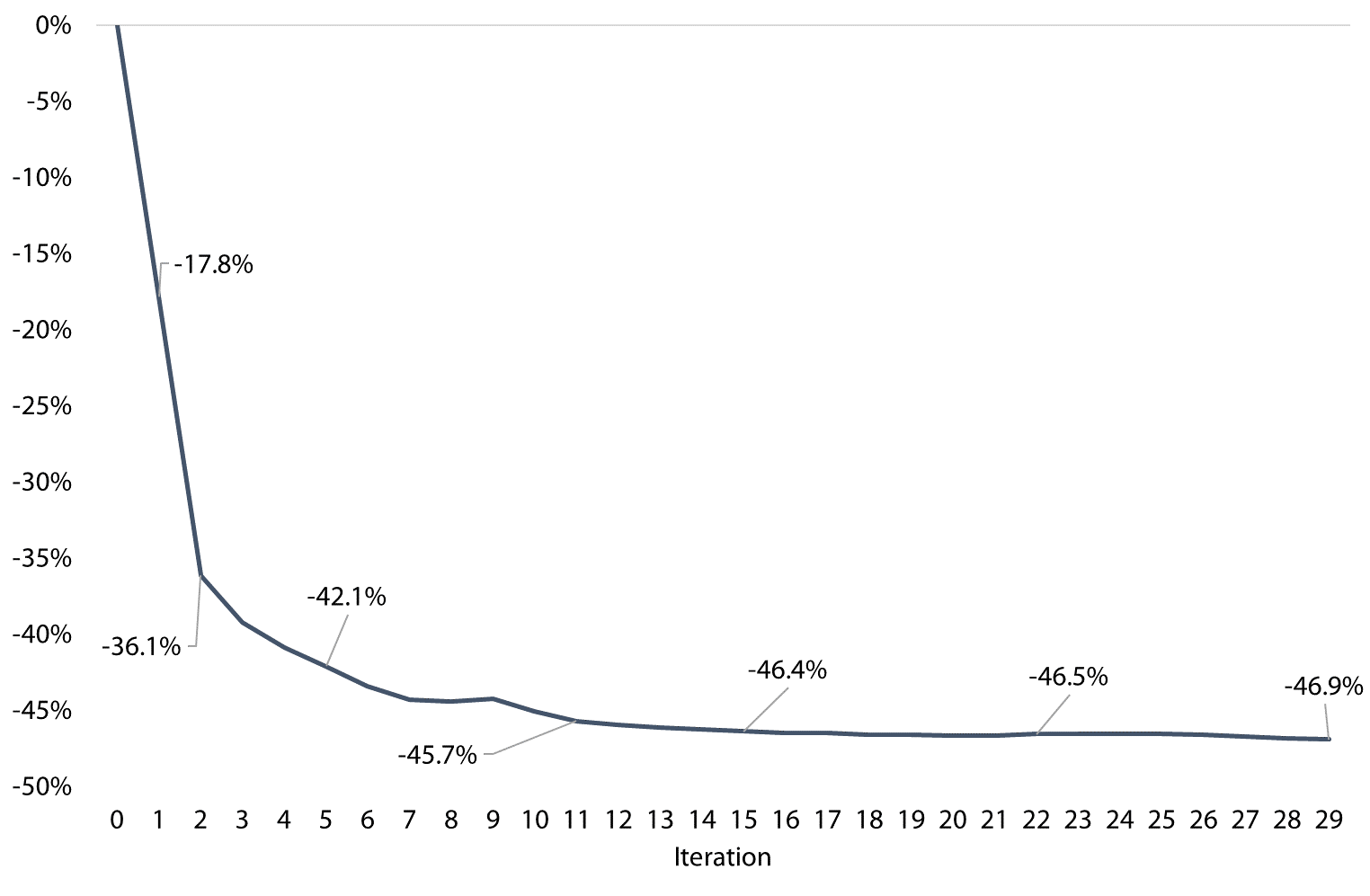}
\caption{Reductions in objective function – Bolivia}
\label{fig: bol_reduction}
\end{figure}

In the first iterations of the method, great reductions on the total system's var requirement are obtained. Fig.~\ref{fig: bol_reduction} illustrates the reductions on total amount invested along the iterative process compared to the solution without progressive hedging (iteration 0) for \textbf{Simulation 12}.

\begin{figure}[!ht]
\centering
\includegraphics[width=5in]{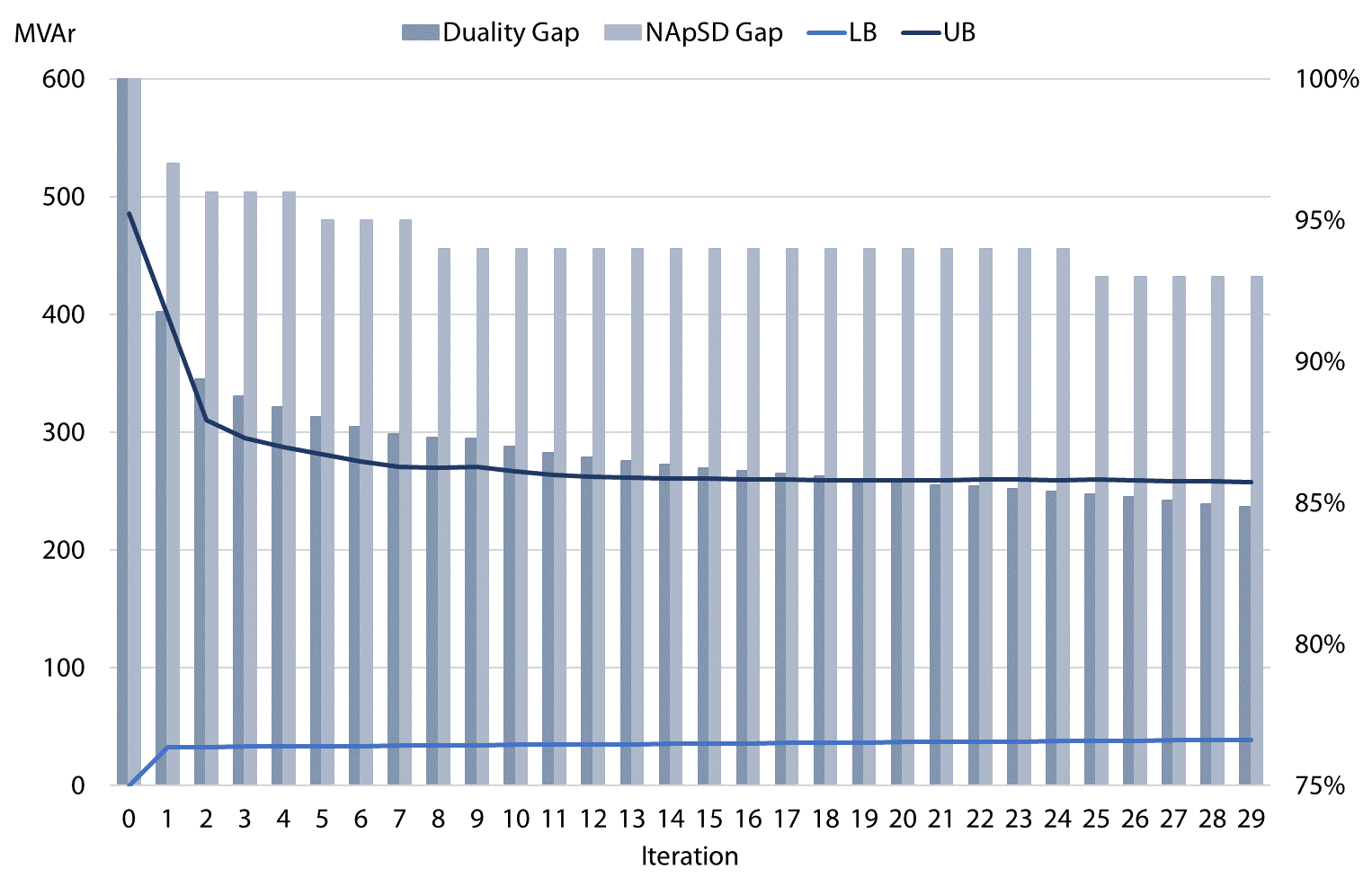}
\caption{LB-UB Convergence scheme – Bolivia}
\label{fig: bol_gap}
\end{figure}

Fig.~\ref{fig: bol_gap} illustrates the convergence of \textbf{Simulation 12}. Note that duality gap is 85\% while the NApSD criterion is 93\%. The approximation between LB and UB are very slow, which wait for convergence would take days of simulation, if it reaches convergence. An approach to accelerate the process aiming the convergence is the parallelization of the simulations or better, the clusterization of the 1,152 operative scenarios considered, reducing the number of scenarios to be solved at each iteration.

\vbox{%
In order to evaluate the impact of the penalty parameter, two more simulations were done with ($\delta_k^{rs}$,$\delta_k^{cs}$) = 0.01. The set to compare now is:

\begin{itemize}
    \item \textbf{Simulation 12:} K = 1;
    \item \textbf{Simulation 14:} K = 0.5;
    \item \textbf{Simulation 15:} K = 2;
\end{itemize}
}

\begin{figure}[!ht]
\centering
\includegraphics[width=5in]{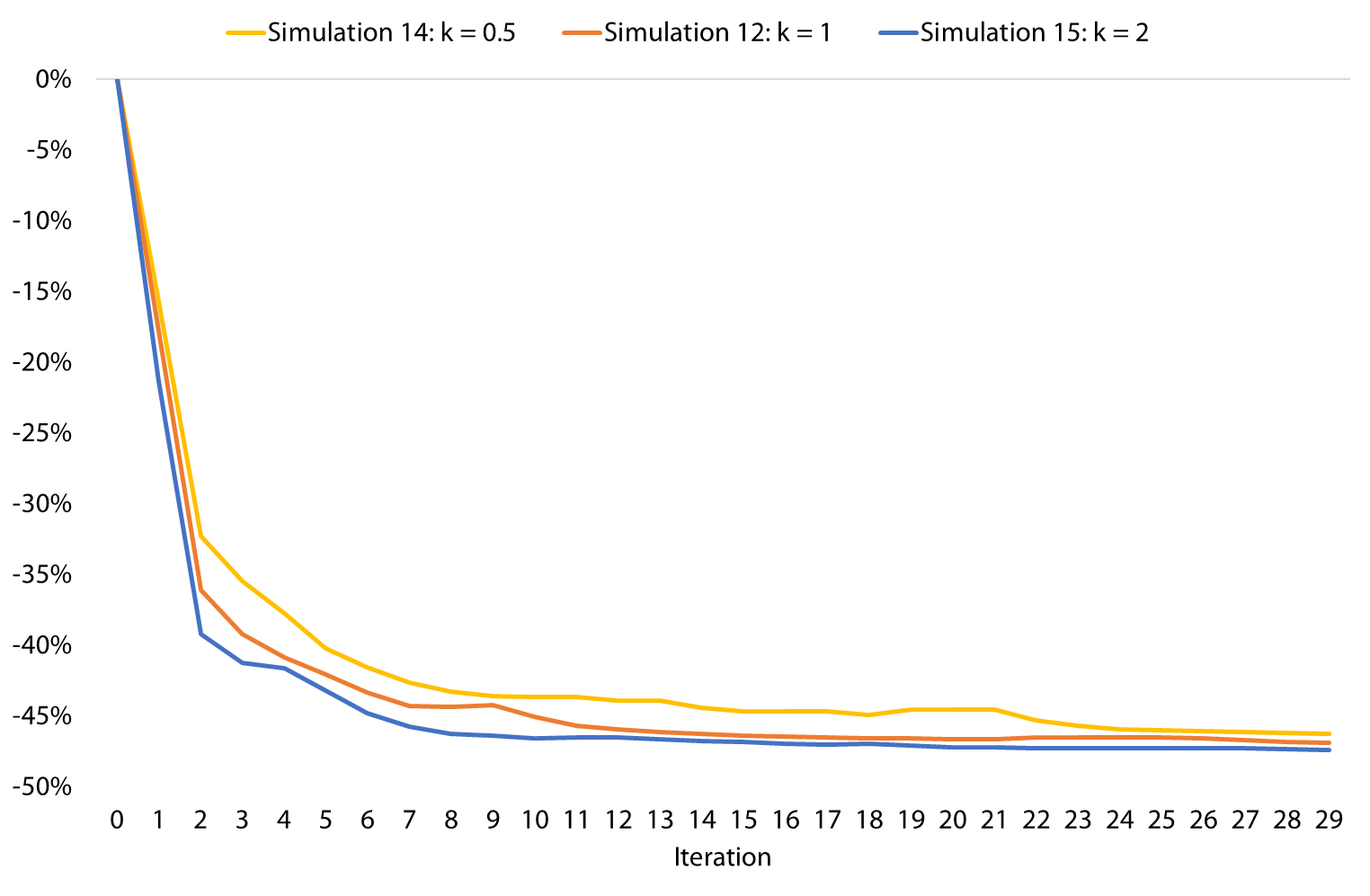}
\caption{Mvar requirement reduction in \% for simulations 12, 14 and 15 – Bolivia}
\label{fig: bol_penaltyanalisys}
\end{figure}

The comparison of the reactive requirement between these three simulations is presented in Fig.~\ref{fig: bol_penaltyanalisys}. It is observed that as the parameter K grows, the more significant is the reduction in the total amount of var required for the first iterations. However, the three simulations seem to stabilize in the same “reduction amount” by the end of the iterative process.

The Table~\ref{tab: CPU Time - Bolivia} summarizes the total CPU time for the whole iterative process from the 5 simulations analyzed for this system. It is observed that all simulations took basically the same computational effort and stopped by the timeout criterium, except \textbf{Simulation 13}, which stopped before 30 iterations because one of the OPFs in the 23rd iteration did not converge. As the number of iterations rise, the CPU time rises as expected, and compared to the simulations with the IEEE 24-Bus, it is observable that the rise in number of scenarios and size of the system makes the total CPU time to rise. 30 iterations for 5 scenarios with the IEEE 24-Bus system took something between 2 and 3 minutes and here takes more than 10 hours. It is also observable that varying the step-size and k seems to have no interference in the CPU Time required to solve the problem with the method’s application.

\begin{table}[!ht]
	\renewcommand{\arraystretch}{1.3}
	\centering
	\caption{Total CPU Time - Bolivian system}
	\label{tab: CPU Time - Bolivia}
	\begin{tabular}{c c c c c c}
		\toprule
			Simulation & t & k & Number of iterations & Total var requirement & CPU Time [min] \\
			\cline{2-6}
            11 & 0.005 & 1 & 30 & 261 & 634 \\
            12 & 0.01 & 1 & 30 & 258 & 638 \\
            13 & 0.1 & 1 & 22 & 486 & 459 \\
            14 & 0.01 & 0.5 & 30 & 261 & 637 \\
            15 & 0.01 & 2 & 30 & 256 & 635 \\
		\bottomrule
	\end{tabular}
\end{table}

\subsection{Colombian Electrical Power System}

\vbox{%
The system consists of:

\begin{itemize}
	\item 1,674 buses;
	\item 1,147 transmission lines;
	\item 1,179 transformers;
	\item 388 generating units;
\end{itemize}
}

\noindent
The list of parameters considered for simulations with the present system is:

\begin{itemize}
	\item Maximum shunt allocation for each bus: 1500 Mvar;
	\item 180 different generation/load set-points considered at each iteration;
	\item Maximum number of iterations: 4;
	\item Convergence gap: 1\%
\end{itemize}

The database for the present simulation came from a real and complete system expansion study. First, a generation expansion plan was obtained for the horizon (2018-2040) \cite{latorre2019stochasticrobust,pereira1991multi,campodonico2003expansion,fern2019stochastic}. Then, it was obtained the transmission expansion plan for the active power transport (reinforcements of transmission lines and transformers). Finally, the reactive power expansion could be then realized after generation and transmission expansions were concluded. For the present simulation, it was considered the analysis for var expansion system for year 2021 (first year of var planning analysis from the study).

The simulation with the Colombian system considered a reduced set of candidates, ($\delta_k^{rs}$,$\delta_k^{cs}$) = 0.01, K = 1 and stopped by timeout criterium (4 iterations in 114 minutes of model's execution). From PH’s application, it is observed that in the first iteration of the method there is a reduction of approximately 27\%, which remains practically constant for the next iterations. Table~\ref{tab: col solution} presents the var requirement at each bus as the convex hulls from the solutions of iterations 0 and 1.

\begin{table}[!ht]
	\renewcommand{\arraystretch}{1.3}
	\centering
	\caption{Solutions for the simulation with the Colombian Eletrical Power System}
	\label{tab: col solution}
	\begin{tabular}{c c c}
		\toprule
			Equipment & Iteration 0 & Iteration 1 \\
			\cline{2-3}
            Reactor at SE Rio Grande   & 10 Mvar &	0 Mvar \\
            Reactor at SE Santa Rosar  & 11.8 Mvar & 11.8 Mvar \\
            Capacitor at SE Mompox     & 10 Mvar & 10 Mvar  \\
            Capacitor at SE Tumaco     & 9.3 Mvar & 9.3 Mvar \\
		\bottomrule
	\end{tabular}
\end{table}

From iteration 0 to iteration 1, it was no longer required the investment in 10 Mvar shunt reactor at SE Rio Grande. This requirement was only in two scenarios (approx. 1\% of the scenarios). Since the PH method induces the variables to approach the reference value set at each iteration, the solution for these two scenarios was the reallocation of the reactor requirement to SE Santa Rosa.

\vbox{%
Considering the real existing modules of shunt equipment (as illustrated in Table~\ref{tab: var cost}), the var expansion plan cost for iteration 0 solution is 17.2 MUS\$ and for iteration 1 is 12.2 MUS\$. The final expansion plan is then:

\begin{itemize}
	\item SE Santa Rosa: 15Mvar Reactor + Connection Bay (5 MUS\%);
	\item SE Mompox: 15Mvar Reactor + Connection Bay (3.6 MUS\%);
	\item SE Tumaco: 10 Mvar Reactor + Connection Bay (3.6 MUS\%).
\end{itemize}
}

Therefore, PH method’s application resulted in big savings (5 MUS\$), reducing the total cost of the var expansion plan by 29\%.

\section{Conclusions}
\label{sec: conclusions}
Currently, due to the high complexity and non-convexity of the optimal nonlinear power flow problem, the analysis of the var expansion planning is done individually for each scenario. In this way, the final investment decision is subject to a superposition of the individual solutions resulting from these scenarios, leading to high expansion costs for the system.

Although it does not guarantee global optimality, the methodology proposed in this work finds a solution that couples investment decisions, referring to each deterministic problem, by imposing a reference value for final decision making. As can be seen from the results, the coupling of the decisions and the proposed methodology resulted in positive results for the var expansion planning problem by reducing the total amount of investments required in the network, meeting all operational requirements.

From simulations with IEEE 24-Bus database, it was concluded, for this particular case, that considering a reduced set of shunt candidates, led to faster convergence and no change in the optimal solution. It was observed that with the introduction of the step-size, the computational effort required doubles, as other optimization problem is required to calculate the lowerbound at each iteration. 

It was also observed from simulations with IEEE 24-Bus and Bolivian system that the chosen parameters can change the way until an optimal solution is achieved (when achieved). With ($\delta_k^{rs}, \delta_k^{rs}$) = 1, as it is usually seen in the literature, the term $w_k^{rs}$ ($q_k^{rs}$ - $q_k^r$) + $w_k^{cs}$ ($q_k^{cs}$ - $q_k^c$) becomes negative in such way that all objective function is negative, and the problem starts to maximize shunt investment instead of minimizing it.

Finally, the simulation from Colombian system presented real gains from the method, where the decision of investing in a 15 Mvar reactor was dismissed, saving 5 MUS\$.

\bibliographystyle{IEEEtran}
\bibliography{myref}

\end{document}